\documentclass[11pt]{article}
\usepackage[numbers,sort&compress]{natbib}
\usepackage{amsfonts}
\usepackage{amsfonts}
\usepackage{amsfonts}
\usepackage{mathrsfs}
\usepackage{bm}
\usepackage[colorlinks,linkcolor=blue,citecolor=blue]{hyperref}
\usepackage{amssymb,amsmath} 
\usepackage{lineno}
 \allowdisplaybreaks

\catcode`!=11
\let\!int\int \def\int{\displaystyle\!int}
\let\!lim\lim \def\lim{\displaystyle\!lim}
\let\!sum\sum \def\sum{\displaystyle\!sum}
\let\!sup\sup \def\sup{\displaystyle\!sup}
\let\!inf\inf \def\inf{\displaystyle\!inf}
\let\!cap\cap \def\cap{\displaystyle\!cap}
\let\!max\max \def\max{\displaystyle\!max}
\let\!min\min \def\min{\displaystyle\!min}
\let\!frac\frac \def\frac{\displaystyle\!frac}
\catcode`!=12

\let\oldsection\section
\renewcommand\section{\setcounter{equation}{0}\oldsection}

\allowdisplaybreaks
\def\pf{\it{Proof.}\rm\quad}

\def\N{\mathbb{N}}

\newtheorem{thm}{Theorem}[section]
\newtheorem{lem}[thm]{Lemma}
\newtheorem{cor}[thm]{Corollary}

\newtheorem{exa}{Example}[section]
\begin{document}
\title {\bf Identities for the multiple zeta (star) values}
\author{
{Ce Xu\thanks{Corresponding author. Email: 15959259051@163.com}}\\[1mm]
\small School of Mathematical Sciences, Xiamen University\\
\small Xiamen
361005, P.R. China}

\date{}
\maketitle \noindent{\bf Abstract } In this paper we prove some new identities for multiple zeta values and multiple zeta star values of arbitrary depth by using the methods of integral computations of logarithm function and iterated integral representations of series. By applying the formulas obtained, we can prove that multiple zeta star values whose indices are the sequences $(\bar 1,\{1\}_m,\bar 1)$ and $(2,\{1\}_m,\bar 1)$ can be expressed polynomially in terms of zeta values, polylogarithms and $\ln(2)$. Finally, we also evaluate several restricted sum formulas involving multiple zeta values.
\\[2mm]
\noindent{\bf Keywords} Multiple zeta value; multiple zeta star value; restricted sum formula.
\\[2mm]
\noindent{\bf AMS Subject Classifications (2010):} 11M06; 11M40; 40B05; 33E20.
\tableofcontents
\section{Introduction}
The study of special values of multiple zeta functions and multiple zeta star functions deals with relations between values at non-zero integer vectors ${\bf s}:=(s_1,\ldots,s_k)$ of sums of the form
\begin{align}\label{1.1}
&\zeta \left( \bf s \right)\equiv\zeta \left( {{s_1}, \ldots ,{s_k}} \right): = \sum\limits_{{n_1} >  \cdots  > {n_k} > 0} {\prod\limits_{j = 1}^k {n_j^{ - \left| {{s_j}} \right|}}{\rm sgn}(s_j)^{n_j}} ,
\end{align}
\begin{align}\label{1.2}
\zeta^\star \left( \bf s \right)\equiv{\zeta ^ \star }\left( {{s_1}, \ldots ,{s_k}} \right): = \sum\limits_{{n_1} \ge  \cdots  \ge {n_k} \ge 1} {\prod\limits_{j = 1}^k {n_j^{ - \left| {{s_j}} \right|}{\rm sgn}(s_j)^{n_j}} } ,
\end{align}
commonly referred to as multiple zeta values and multiple zeta star values \cite{BBBL1997,BBBL2001,BB2003,O1999,Y2009}, respectively, where
\[{\mathop{\rm sgn}} \left( {{s_j}} \right): = \left\{ {\begin{array}{*{20}{c}}
   {1,} & {{s_j} > 0,}  \\
   { - 1,} & {{s_j} < 0.}  \\
\end{array}} \right.\]
 We are primarily interested in positive integer values of the arguments $s_1,\ldots,s_k$, in which case it is seen that the condition $s_1>1$ is necessary for both sums (\ref{1.1}) and (\ref{1.2}) to converge. Of course, if ${\sigma_1} =-1$, then we can allow $s_1=1$. Throughout the paper we will use $\bar p$ to denote a negative entry $s_j=-p$. For example,
 \[\zeta \left( {{{\bar s}_1},{s_2}} \right) = \zeta \left( { - {s_1},{s_2}} \right),\ {\zeta ^ \star }\left( {{{\bar s}_1},{s_2},{{\bar s}_3}} \right) = {\zeta ^ \star }\left( { - {s_1},{s_2}, - {s_3}} \right).\]
We call $l({\bf s}):=k$ the depth of (\ref{1.1}), (\ref{1.2}) and $\left| {\bf s} \right|: = \sum\limits_{j = 1}^k {\left| {{s_j}} \right|} $ the weight. For convenience we set $\zeta \left( \emptyset  \right) = {\zeta ^ \star }\left( \emptyset  \right) = 1$ and ${\left\{ {{s_1}, \ldots ,{s_j}} \right\}_d}$ to be the set formed by repeating the composition $\left( {{s_1}, \ldots ,{s_j}} \right)$ $d$ times.

In this paper we let $\zeta_n \left( {{s_1}, \ldots ,{s_k}} \right)$ and $\zeta_n^\star \left( {{s_1}, \ldots ,{s_k}} \right)$ to denote the multiple harmonic number (also called the partial sums of multiple zeta value) and multiple harmonic star number (also called the partial sums of multiple zeta star value), respectively, which are defined by
\begin{align}\label{a1}
{\zeta _n}\left( {{s_1},{s_2}, \cdots ,{s_k}} \right): = \sum\limits_{n \ge {n_1} > {n_2} >  \cdots  > {n_k} \ge 1} {\prod\limits_{j = 1}^k {n_j^{ - \left| {{s_j}} \right|}} {\rm{sgn}}{{({s_j})}^{{n_j}}}} ,
\end{align}
\begin{align}\label{a2}
\zeta _n^ \star \left( {{s_1},{s_2}, \cdots ,{s_k}} \right): = \sum\limits_{n \ge {n_1} \ge {n_2} \ge  \cdots  \ge {n_k} \ge 1} {\prod\limits_{j = 1}^k {n_j^{ - \left| {{s_j}} \right|}} {\rm{sgn}}{{({s_j})}^{{n_j}}}} ,
\end{align}
when $n<k$, then ${\zeta _n}\left( {{s_1},{s_2}, \cdots ,{s_k}} \right)=0$, and ${\zeta _n}\left(\emptyset \right)={\zeta^\star _n}\left(\emptyset \right)=1$.

A good deal of work on multiple zeta (star) values has focused on the problem of determining when `complicated' sums can be expressed in terms of `simpler' sums. Thus, researchers are interested in determining which sums can be expressed in terms of other sums of smaller depth.
The origin of multiple zeta (star) values goes back to the correspondence of Euler with Goldbach in 1742-1743 (see \cite{H2007}). Euler studied double zeta values and established some important relation formulas for them. For example, he proved that
\begin{align}\label{1.3}
2{\zeta ^ \star }\left( {k,1} \right) = \left( {k + 2} \right)\zeta \left( {k + 1} \right) - \sum\limits_{i = 1}^{k - 2} {\zeta \left( {k - i} \right)\zeta \left( {i + 1} \right)},\ k\geq2,
\end{align}
which, in particular, implies the simplest but nontrivial relation:
${\zeta ^ \star }\left( {2,1} \right)\; = 2\zeta \left( 3 \right)$ or equivalently $\zeta \left( {2,1} \right)\; = \zeta \left( 3 \right)$.
Moreover, he conjectured that the double zeta values would be reducible whenever weight is odd, and even gave what he hoped to be the general formula. The conjecture was first proved by Borweins and Girgensohn \cite{BBG1995}. For some interesting results on generalized double zeta values (also called Euler sums), see \cite{BBG1994,FS1998}.

The systematic study of multiple zeta values began in the early 1990s with the works of Hoffman \cite{H1992}, Zagier \cite{DZ1994} and Borwein-Bradley-Broadhurst \cite{BBBL1997} and has continued with increasing attention in recent years (see \cite{CCE2016,CML2016,E2009,E2016}). The first systematic study of reductions up to depth 3 was carried out by Borwein and Girgensohn \cite{BG1996}, where the authors proved that if $p+q+r$ is even or less than or equal to 10 or $p+q+r=12$, then triple zeta values $\zeta \left( {q,p,r} \right)$ (or $\zeta^\star \left( {q,p,r} \right)$) can be expressed as a rational linear combination of products of zeta values and double zeta values. The set of the the multiple zeta values has a rich algebraic structure given by the shuffle and the stuffle relations \cite{H1992,M2009,DZ2012}.

Additionally, it has been discovered in the course of the years that many multiple zeta (star) values admit expressions involving finitely many zeta values and polylogarithms, that are to say values of the Riemann zeta function and polylogarithm function,
\begin{align*}
&\zeta \left( s \right): = \sum\limits_{n = 1}^\infty  {\frac{1}{{{n^s}}}} \;\left( {{\mathop{\Re}\nolimits} \left( s \right) > 1} \right),\\
&{\rm Li}_s\left( x \right): = \sum\limits_{n = 1}^\infty  {\frac{{{x^n}}}{{{n^s}}}} \;\left( {{\mathop{\Re}\nolimits} \left( s \right) \ge 1,\;x \in \left[ { - 1,1} \right)} \right)
\end{align*}
at positive integer arguments and $x=1/2$. The relations between multiple zeta (star) values and the values of the Riemann zeta function and polylogarithm function have been studied by many authors, see \cite{BBBL1997,BBBL2001,BB2003,L2012,KO2010,O1999,KP2013,KPZ2016,X2017,Y2009,DZ2012}. For example, Zagier \cite{DZ2012} proved that the multiple zeta star values ${\zeta ^ \star \left( {{{\left\{ 2 \right\}}_a},3,{{\left\{ 2 \right\}}_b}} \right)} $ and multiple zeta values ${\zeta \left( {{{\left\{ 2 \right\}}_a},3,{{\left\{ 2 \right\}}_b}} \right)} $, where $a,b\in \N_0:=\N\cup \{0\}=\{0,1,2,\ldots\}$, are reducible to polynomials in zeta values, and gave explicit formulae. Hessami Pilehrood et al. \cite{KP2013} and Li \cite{L2012} provide two new proofs of Zagier's formula for ${\zeta ^ \star \left( {{{\left\{ 2 \right\}}_a},3,{{\left\{ 2 \right\}}_b}} \right)} $ based on a finite identity for partial sums of the zeta-star series and hypergeometric series computations, respectively.

The purpose of the present paper is to establish some new family of identities for multiple integral representation of series. Then, we apply it to obtain a family of identities relating multiple zeta (star) values to zeta valus and polylogarithms. Specially, we present some new recurrence relations for multiple zeta star values whose indices are the sequences $(\bar 1,\{1\}_m,\bar 1)$, $(2,\{1\}_m,\bar 1)$, and prove that the multiple zeta values $\zeta \left( {\bar 1,{{\left\{ 1 \right\}}_{m}},\bar 1} \right)$ and $\zeta \left( {\bar 1,{{\left\{ 1 \right\}}_{m}},\bar 1,\bar 1,{{\left\{ 1 \right\}}_{k}}} \right)$
can be expressed as a rational linear combination
of zeta values, polylogarithms and $\ln(2)$. Moreover, we also evaluate several restricted sum formulas involving multiple zeta values.
\section{Evaluation of the multiple zeta star values $\zeta^\star(\bar 1,\{1\}_m,\bar 1)$ and $\zeta^\star(2,\{1\}_m,\bar 1)$}
In this section, we will show that the alternating multiple zeta star values $\zeta^\star(\bar 1,\{1\}_m,\bar 1)$ and $\zeta^\star(2,\{1\}_m,\bar 1)$ satisfy certain recurrence relations that allow us to write them in terms of zeta values, polylogarithms and $\ln(2)$.
\subsection{Four lemmas and proofs}
Now, we need four lemmas which will be useful in proving Theorems 2.4 and 2.5.
\begin{lem}(\cite{Xu2017}) \label{lem2.1}
For $n\in \N,\ m\in \N_0$ and $x\in[-1,1)$, the following relation holds:
\begin{align}\label{2.1}
\int\limits_0^x {{t^{n - 1}}{{\ln }^m}\left( {1 - t} \right)} dt&=\frac{1}{n}\left( {{x^n} - 1} \right){\ln ^m}\left( {1 - x} \right)+ m!\frac{{{{\left( { - 1} \right)}^m}}}{n}\sum\limits_{1 \le {k_m} \le  \cdots  \le {k_1} \le n} {\frac{{{1}}}{{{k_1} \cdots {k_m}}}}\nonumber \\
&\quad - \frac{1}{n}\sum\limits_{i = 1}^{m } {{{\left( { - 1} \right)}^{i - 1}}i!\left( {\begin{array}{*{20}{c}}
   m  \\
   i  \\
\end{array}} \right){{\ln }^{m - i}}\left( {1 - x} \right)} \sum\limits_{1 \le {k_i} \le  \cdots  \le {k_1} \le n} {\frac{{{x^{{k_i}}} - 1}}{{{k_1} \cdots {k_i}}}},
\end{align}
where if $m=0$, then
\[\sum\limits_{_{1 \le {k_m} \le  \cdots  \le {k_1} \le n}} {\frac{{f\left( {{k_m}} \right)}}{{{k_1} \cdots {k_m}}}}  = f\left( n \right).\]

\end{lem}
\pf The proof is by induction on $m$.
Define $J\left( {n,m;x} \right): = \int\limits_0^x {{t^{n - 1}}{{\ln }^m}\left( {1 - t} \right)} dt$, for $m=0$, by a simple calculation, we can arrive at the conclusion that
\[J\left( {n,0;x} \right) = \int\limits_0^x {{t^{n - 1}}dt}  = \frac{{{x^n}}}{n}\]
and the formula is true. For $m\geq 1$ we proceed as follows.
By using integration by parts we have the following recurrence relation
\begin{align}\label{2.2}
J\left( {n,m;x} \right) = \frac{1}{n}\left( {{x^n} - 1} \right){\ln ^m}\left( {1 - x} \right) - \frac{m}{n}\sum\limits_{k = 1}^n {J\left( {k,m - 1;x} \right)}.
\end{align}
Suppose the lemma holds for $m-1$. Then the inductive hypothesis implies that the integral (\ref{2.2}) is equal to
\begin{align}\label{2.3}
J\left( {n,m;x} \right) =& \frac{1}{n}\left( {{x^n} - 1} \right) {\ln ^m}\left( {1 - x} \right)+ m!\frac{{{{\left( { - 1} \right)}^m}}}{n}\sum\limits_{j = 1}^n \frac 1{j}{\sum\limits_{1 \le {k_{m - 1}} \le  \cdots  \le {k_1} \le j} {\frac{{{1}}}{{{k_1} \cdots {k_{m - 1}}}}} }\nonumber \\
& + \frac{m}{n}\sum\limits_{i = 1}^{m - 1} {{{\left( { - 1} \right)}^{i - 1}}i!\left( {\begin{array}{*{20}{c}}
   {m - 1}  \\
   i  \\
\end{array}} \right){{\ln }^{m - 1 - i}}\left( {1 - x} \right)}\sum\limits_{j = 1}^n \frac 1{j} \sum\limits_{1 \le {k_i} \le  \cdots  \le {k_1} \le j} {\frac{{{x^{{k_i}}} - 1}}{{{k_1} \cdots {k_i}}}}\nonumber \\
& - \frac{m}{n}{\ln ^{m - 1}}\left( {1 - x} \right)\sum\limits_{j = 1}^n {\frac{{{x^j} - 1}}{j}} .
\end{align}
Thus, by a direct calculation, we may deduce the desired result. This completes the proof of Lemma 2.1. \hfill$\square$

\begin{lem}\label{blem1}
 For integers $m\geq 0, n\geq 0$ and $x\in (0,1)$, then the following integral identity holds:
\begin{align}\label{b1}
\int\limits_0^x {{t^{n}}{{\left( {\ln t} \right)}^m}} dt = \sum\limits_{l = 0}^m {l!\left( {\begin{array}{*{20}{c}}
   m  \\
   l  \\
\end{array}} \right)\frac{{{{\left( { - 1} \right)}^l}}}{{{(n+1)^{l + 1}}}}{x^{n+1}}{{{\ln^{m - l} (x)} }}}.
\end{align}
\end{lem}
\pf This follows by induction on $m$, using the formula
\[\int\limits_0^{{x}} {{t^{n }}{{\ln }^m}(t)} dt = \frac{{{x^{n+1}}}}{n+1}{\ln ^m}(x) - \frac{m}{n+1}\int\limits_0^{{x}} {{t^{n}}{{\ln }^{m - 1}}(t)} dt,\]
that comes from integration by parts.\hfill$\square$

\begin{lem} (\cite{X2017}) \label{lem2.2}
For positive integer $m$, the following identity holds:
\begin{align}\label{2.4}
\int\limits_0^1 {\frac{{{{\ln }^m}\left( {1 + t} \right)\ln \left( {1 - t} \right)}}{{1 + t}}} dt= &\frac{1}{{m + 1}}{\ln ^{m + 2}}(2) - \zeta \left( 2 \right){\ln ^m}(2)\nonumber\\
& - \sum\limits_{k = 1}^m {\left( {\begin{array}{*{20}{c}}
   m  \\
   k  \\
\end{array}} \right){{\left( { - 1} \right)}^{k + 1}}\left\{ \begin{array}{l}
 \sum\limits_{l = 1}^k {l!\left( {\begin{array}{*{20}{c}}
   k  \\
   l  \\
\end{array}} \right){{\left( {\ln^{m - l} (2)} \right)}}{\rm{L}}{{\rm{i}}_{l + 2}}\left( {\frac{1}{2}} \right)}  \\
  - k!{\left( {\ln^{m - k} (2)} \right)}\zeta \left( {k + 2} \right) \\
 \end{array} \right\}}.
\end{align}
\end{lem}
\pf This proof is based on the simple variable substitution $1+t=2x$ and the help of formula (\ref{b1}). We have
\begin{align*}
 \int\limits_0^1 {\frac{{{{\ln }^m}\left( {1 + t} \right)\ln \left( {1 - t} \right)}}{{1 + t}}} dt\rm{ = }&\int\limits_{1/2}^1 {\frac{{{{\ln }^m}\left( {2x} \right)\ln \left( {2 - 2x} \right)}}{x}dx}  \\
  =& \sum\limits_{k = 0}^m {\left( {\begin{array}{*{20}{c}}
   m  \\
   k  \\
\end{array}} \right){{\ln }^{m - k + 1}}\left( 2 \right)\int\limits_{1/2}^1 {\frac{{{{\ln }^k}\left( x \right)}}{x}dx} }  \\
 & + \sum\limits_{k = 0}^m {\left( {\begin{array}{*{20}{c}}
   m  \\
   k  \\
\end{array}} \right){{\ln }^{m - k}}\left( 2 \right)\int\limits_{1/2}^1 {\frac{{{{\ln }^k}\left( x \right)\ln \left( {1 - x} \right)}}{x}dx} } .
\end{align*}
Then by using identity (\ref{b1}) and the power series expansion of logarithm function
\[\ln \left( {1 - x} \right) =  - \sum\limits_{n = 1}^\infty  {\frac{{{x^n}}}{n}} ,\quad x\in [-1,1)\]
we easily obtain the formula (\ref{2.4}), which completes the proof of Lemma \ref{lem2.2}.\hfill$\square$

\begin{lem} (\cite{X2017}) \label{lem2.3}
For integer $m\geq1$ and $x \in \left[ {0,1} \right]$, the following identity holds:
\begin{align}\label{2.5}
\int\limits_0^x {\frac{{{{\ln }^m}\left( {1 + t} \right)}}{t}} dt = &\frac{1}{{m + 1}}{\ln ^{m + 1}}\left( {1 + x} \right) + m!\left( {\zeta \left( {m + 1} \right) - {\rm Li}{_{m + 1}}\left( {\frac{1}{{1 + x}}} \right)} \right)\nonumber \\&- m!\sum\limits_{j = 1}^m {\frac{{{{\ln }^{m - j + 1}}\left( {1 + x} \right)}}{{\left( {m - j + 1} \right)!}}} {\rm Li}{_j}\left( {\frac{1}{{1 + x}}} \right).
\end{align}
\end{lem}
\pf We note that applying the change of variable $t\rightarrow w-1$ to the left hand side of (\ref{2.5}), which can be rewritten as
\begin{align*}
\int\limits_0^x {\frac{{{{\ln }^m}\left( {1 + t} \right)}}{t}} dt\mathop  = \limits^{w = 1 +t}  & \int\limits_1^{1 + x} {\frac{{{{\ln }^m}(w)}}{{w - 1}}} dw\mathop  = \limits^{u = {w^{ - 1}}} {\left( { - 1} \right)^{m + 1}}\int\limits_1^{{{(1 + x)^{-1}}}} {\frac{{{{\ln }^m}(u)}}{{u - {u^2}}}} du\\
           =& {\left( { - 1} \right)^{m + 1}}\left\{ {\int\limits_1^{{{(1 + x)^{-1}}}} {\frac{{{{\ln }^m}(u)}}{u}} du + \int\limits_1^{{{(1 + x)^{-1}}}} {\frac{{{{\ln }^m}(u)}}{{1 - u}}} du} \right\}\\
           =&\frac{1}{{m + 1}}{\ln ^{m + 1}}\left( {1 + x} \right) + {\left( { - 1} \right)^{m + 1}}\int\limits_1^{{{(1 + x)^{-1}}}} {\frac{{{{\ln }^m}(u)}}{{1 - u}}} du .
\end{align*}
Then with the help of formula (\ref{b1}) we may deduce the evaluation (\ref{2.5}).\hfill$\square$

\subsection{Two theorems and proofs}
Now, we state our main results, the main results in this section are the following two theorems.
\begin{thm}\label{thm1}
For positive integer $m$, we have the recurrence relation
\begin{align}\label{2.6}
\zeta^\star \left( {\bar 1,{{\left\{ 1 \right\}}_{m }},\bar 1} \right)= &\frac{{{{\left( { - 1} \right)}^{m }}}}{{ {m }!}}\zeta \left( 2 \right){\ln ^{m}}(2)\nonumber\\
& - \frac{{{{\left( { - 1} \right)}^{m }}}}{{(m+1)!}}\sum\limits_{i = 1}^{m} {{{\left( { - 1} \right)}^{i + 1}}i!\left( {\begin{array}{*{20}{c}}
   m+1  \\
   i  \\
\end{array}} \right){{\left( {\ln 2} \right)}^{m +1- i}}} \left\{ {\zeta^\star \left( {\bar 1,{{\left\{ 1 \right\}}_{i - 1}},\bar 1} \right) - \zeta^\star \left( {\bar 1,{{\left\{ 1 \right\}}_i}} \right)} \right\}\nonumber\\
&+ \frac{{{{\left( { - 1} \right)}^{m}}}}{{{m}!}}\sum\limits_{k = 1}^{m} {\left( {\begin{array}{*{20}{c}}
   {m}  \\
   k  \\
\end{array}} \right){{\left( { - 1} \right)}^{k + 1}}\left\{ \begin{array}{l}
 \sum\limits_{l = 1}^k {l!\left( {\begin{array}{*{20}{c}}
   k  \\
   l  \\
\end{array}} \right){{\left( {\ln 2} \right)}^{m - l}}{\rm{L}}{{\rm{i}}_{l + 2}}\left( {\frac{1}{2}} \right)}  \\
  - k!{\left( {\ln 2} \right)^{m- k}}\zeta \left( {k + 2} \right) \\
 \end{array} \right\}},
\end{align}
where \[{\zeta ^ \star }\left( {\bar 1,\bar 1} \right) = \frac{{\zeta \left( 2 \right) + {{\ln }^2}\left( 2 \right)}}{2}.\]
\end{thm}
\pf Multiplying (\ref{2.1}) by ${\left( { - 1} \right)^{n - 1}}$ and summing with respect to $n$, the result is
\begin{align}\label{2.7}
&\int\limits_0^x {\frac{{{{\ln }^m}\left( {1 - t} \right)}}{{1 + t}}} dt = \sum\limits_{n = 1}^\infty  {{{\left( { - 1} \right)}^{n - 1}}\int\limits_0^x {{t^{n - 1}}{{\ln }^m}\left( {1 - t} \right)} dt}\nonumber \\
&= {\ln ^m}\left( {1 - x} \right)\sum\limits_{n = 1}^\infty  {\frac{{{x^n} - 1}}{n}{{\left( { - 1} \right)}^{n - 1}}}  + m!{\left( { - 1} \right)^m}\sum\limits_{n = 1}^\infty  {\frac{{{{\left( { - 1} \right)}^{n - 1}}}}{n}\sum\limits_{1 \le {k_m} \le  \cdots  \le {k_1} \le n} {\frac{{{x^{{k_m}}}}}{{{k_1} \cdots {k_m}}}} }\nonumber \\
& \quad- \sum\limits_{i = 1}^{m - 1} {{{\left( { - 1} \right)}^{i - 1}}i!\left( {\begin{array}{*{20}{c}}
   m  \\
   i  \\
\end{array}} \right){{\ln }^{m - i}}\left( {1 - x} \right)} \sum\limits_{n = 1}^\infty  {\frac{{{{\left( { - 1} \right)}^{n - 1}}}}{n}\sum\limits_{1 \le {k_i} \le  \cdots  \le {k_1} \le n} {\frac{{{x^{{k_i}}} - 1}}{{{k_1} \cdots {k_i}}}} }.
\end{align}
On the other hand, by integration by parts, we obtain the formula
\begin{align*}
&\mathop {\lim }\limits_{x \to  - 1} \left\{ {\int\limits_0^x {\frac{{{{\ln }^m}\left( {1 - t} \right)}}{{1 + t}}} dt - {{\ln }^m}\left( {1 - x} \right)\ln \left( {1 + x} \right)} \right\}\\
& = \mathop {\lim }\limits_{x \to  - 1} \left\{ {m\int\limits_0^x {\frac{{{{\ln }^{m - 1}}\left( {1 - t} \right)\ln \left( {1 + t} \right)}}{{1 - t}}} dt} \right\}\\
& = m\int\limits_0^{ - 1} {\frac{{{{\ln }^{m - 1}}\left( {1 - t} \right)\ln \left( {1 + t} \right)}}{{1 - t}}} dt\\
& =  - m\int\limits_0^1 {\frac{{{{\ln }^{m - 1}}\left( {1 + t} \right)\ln \left( {1 - t} \right)}}{{1 + t}}} dt,
\end{align*}
Hence, replacing $m-1$ by $m$ and letting $x$ approach $-1$ in (\ref{2.7}), and combining it with (\ref{2.4}), we deduce (\ref{2.6}). \hfill$\square$\\
Noting that, taking $x\rightarrow 1$ in (\ref{2.7}), we have the result
\[\int\limits_0^1 {\frac{{{{\ln }^m}\left( {1 - t} \right)}}{{1 + t}}} dt = m!{\left( { - 1} \right)^{m - 1}}\zeta^\star \left( {\bar 1,{{\left\{ 1 \right\}}_m}} \right).\]
On the other hand, applying the change of variable $t\rightarrow 1-x$ to integral above and using (\ref{b1}), we obtain
\[\int\limits_0^1 {\frac{{{{\ln }^m}\left( {1 - t} \right)}}{{1 + t}}} dt = {\left( { - 1} \right)^m}m!{\rm{L}}{{\rm{i}}_{m + 1}}\left( {\frac{1}{2}} \right).\]
Therefore, we conclude that

\begin{align}\label{2.8}
\zeta^\star \left( {\bar 1,{{\left\{ 1 \right\}}_m}} \right) =  - {\rm{L}}{{\rm{i}}_{m + 1}}\left( {\frac{1}{2}} \right),\quad m\in\N_0.
\end{align}

\begin{thm}\label{thm2}
For positive integer $m$, we have the recurrence relation
\begin{align}\label{2.9}
\zeta^\star \left( {2,{{\left\{ 1 \right\}}_{m}},\bar 1} \right) =& \frac{{m + 2}}{{\left( {m + 3} \right)!}}{\left( { - 1} \right)^{m }}{\ln ^{m + 3}}(2) + \left( {m + 2} \right){\left( { - 1} \right)^{m }}\left( {\zeta \left( {m + 3} \right) - {\rm{L}}{{\rm{i}}_{m + 3}}\left( {\frac{1}{2}} \right)} \right)\nonumber\\
& - \left( {m + 2} \right){\left( { - 1} \right)^{m }}\sum\limits_{j = 1}^{m + 2} {\frac{{{{ {\ln^{m + 3 - j} (2)}}}}}{{\left( {m + 3 - j} \right)!}}} {\rm{L}}{{\rm{i}}_j}\left( {\frac{1}{2}} \right) - \frac{3}{2}\frac{{{{\left( { - 1} \right)}^{m}}}}{{(m+1)!}}\zeta \left( 2 \right){ {\ln^{m+1} (2) }}\nonumber\\
& - \frac{{{{\left( { - 1} \right)}^{m}}}}{{(m+1)!}}\sum\limits_{i = 1}^{m} {{{\left( { - 1} \right)}^{i - 1}}i!\left( {\begin{array}{*{20}{c}}
   m+1  \\
   i  \\
\end{array}} \right){{\left( {\ln^{m+1 - i}(2)} \right)}}}\nonumber\\
&\quad\quad\quad\quad\quad\quad\times  \left\{ {\zeta^\star \left( {2,{{\left\{ 1 \right\}}_{i - 1}},\bar 1} \right) - \zeta^\star \left( {2,{{\left\{ 1 \right\}}_i}} \right)} \right\},
\end{align}
where \[{\zeta ^ \star }\left( {2,\bar 1} \right) = \frac{1}{4}\zeta \left( 3 \right) - \frac{3}{2}\zeta \left( 2 \right)\ln \left( 2 \right).\]
\end{thm}
\pf Similarly as in the proof of Theorem \ref{thm1}, we consider the integral
\begin{align}\label{2.10}
\int\limits_0^{ - 1} {\frac{{{{\ln }^{m + 1}}\left( {1 - t} \right)}}{t}} dt =  - \sum\limits_{n = 1}^\infty  {\frac{1}{n}\int\limits_0^{ - 1} {{t^{n - 1}}{{\ln }^m}\left( {1 - t} \right)dt} }.
\end{align}
Putting $x\rightarrow -1$ in (\ref{2.1}), we deduce that
\begin{align}\label{2.11}
\int_0^{ - 1} {{t^{n - 1}}\ln^m \left( {1 - t} \right)} dt =& \frac{1}{n}({\ln ^m}(2))\left( {{{\left( { - 1} \right)}^n} - 1} \right) + m!\frac{{{{\left( { - 1} \right)}^m}}}{n}{\zeta^\star _n}\left( {{{\left\{ 1 \right\}}_{m-1}},\bar 1} \right)\nonumber\\
& - \frac{1}{n}\sum\limits_{i = 1}^{m - 1} {{{\left( { - 1} \right)}^{i - 1}}i!\left( {\begin{array}{*{20}{c}}
   m  \\
   i  \\
\end{array}} \right)} {(\ln ^{m - i}}(2))\left\{ {{\zeta^\star _n}\left( {{{\left\{ 1 \right\}}_{i - 1}},\bar 1} \right) - {\zeta^\star _n}\left( {{{\left\{ 1 \right\}}_i}} \right)} \right\}.
\end{align}
Setting $x\rightarrow 1$ in (\ref{2.5}) and combining (\ref{2.10}) with (\ref{2.11}), the replacing $m-1$ by $m$, we obtain the desired result (\ref{2.9}).\hfill$\square$\\
By considering the case $m=2$ in (\ref{2.6}) and (\ref{2.9}) we get
\[\zeta^\star \left( {\bar 1,1,\bar 1} \right){\rm{ = }}\frac{1}{8}\zeta \left( 3 \right){\rm{ + }}\frac{1}{2}\zeta \left( 2 \right)\ln(2) - \frac{1}{6}{\ln ^3}(2),\]
\[\zeta^\star \left( {2,1,\bar 1} \right) = \frac{1}{8}{\ln ^4}(2) + 3{\rm{L}}{{\rm{i}}_4}\left( {\frac{1}{2}} \right) - 3\zeta \left( 4 \right) - \frac{3}{2}\zeta \left( 2 \right){\ln ^2}(2) + \frac{7}{8}\zeta \left( 3 \right)\ln (2).\]
From \cite{Xu2017} we have the result
\[\zeta^\star \left( {2,{{\left\{ 1 \right\}}_m}} \right) = \left( {m + 1} \right)\zeta \left( {m + 2} \right),\quad m\in\N_0.\]
Therefore, from Theorem \ref{thm1} and Theorem \ref{thm2}, we know that the alternating multiple zeta star values $\zeta^\star(\bar 1,\{1\}_m,\bar 1)$ and $\zeta^\star(2,\{1\}_m,\bar 1)$ can be expressed as a rational linear combination of zeta values, polylogarithms and $\ln(2)$.

\section{Some results on multiple zeta values}
In this section, we use certain multiple integral representations to evaluate several multiple zeta values. We need the following lemma.
\begin{lem} (\cite{X2017,Xu2017})\label{lem3.1}
For integer $k>0$ and $x\in [-1,1)$, we have that
\begin{align}\label{3.1}
{\ln ^k}\left( {1 - x} \right) = {\left( { - 1} \right)^k}k!\sum\limits_{n = 1}^\infty  {\frac{{{x^n}}}{n}{\zeta _{n - 1}}\left( {{{\left\{ 1 \right\}}_{k - 1}}} \right)},
\end{align}
\begin{align}\label{3.2}
s\left( {n,k} \right) = \left( {n - 1} \right)!{\zeta _{n - 1}}\left( {{{\left\{ 1 \right\}}_{k - 1}}} \right).
\end{align}
where ${s\left( {n,k} \right)}$ is called (unsigned) Stirling number of the first kind (see \cite{L1974}).
The Stirling numbers ${s\left( {n,k} \right)}$ of the first kind satisfy a recurrence relation in the form
\[s\left( {n,k} \right) = s\left( {n - 1,k - 1} \right) + \left( {n - 1} \right)s\left( {n - 1,k} \right),\;\;n,k \in \N,\]
with $s\left( {n,k} \right) = 0,n < k,s\left( {n,0} \right) = s\left( {0,k} \right) = 0,s\left( {0,0} \right) = 1$.
\end{lem}
\pf  The proof is based on the two identities
\begin{align}\label{c1}
{\ln ^{k{\rm{ + }}1}}\left( {1 - x} \right){\rm{ = }} - \left( {k + 1} \right)\int\limits_0^x {\frac{{{{\ln }^k}\left( {1 - t} \right)}}{{1 - t}}dt},\quad k\in \N_0
\end{align}
and
\begin{align}\label{c2}
{\ln ^k}\left( {1 - x} \right) = {\left( { - 1} \right)^k}k!\sum\limits_{n = k}^\infty  {\frac{{s\left( {n,k} \right)}}{{n!}}{x^n}} ,\: - 1 \le x < 1.
\end{align}
Applying the induction hypothesis and Cauchy product formula, we arrive at
\begin{align*}
{\ln ^{k{\rm{ + }}1}}\left( {1 - x} \right){\rm{ = }}& - \left( {k + 1} \right)\int\limits_0^x {\frac{{{{\ln }^k}\left( {1 - t} \right)}}{{1 - t}}dt} \\
& = {\left( { - 1} \right)^{k + 1}}\left( {k + 1} \right)!\sum\limits_{n = 1}^\infty  {\frac{1}{{n + 1}}\sum\limits_{i = 1}^n {\frac{{{\zeta _{i - 1}}\left( {{{\left\{ 1 \right\}}_{k - 1}}} \right)}}{i}} } {x^{n + 1}}\\
& = {\left( { - 1} \right)^{k + 1}}\left( {k + 1} \right)!\sum\limits_{n = 1}^\infty  {\frac{{{\zeta _n}\left( {{{\left\{ 1 \right\}}_k}} \right)}}{{n + 1}}} {x^{n + 1}}.
\end{align*}
Nothing that ${\zeta _n}\left( {{{\left\{ 1 \right\}}_k}} \right) = 0$ when $n<k$. Thus, we can deduce (\ref{3.1}). Then, comparing the coefficients of $x^n$ in (\ref{3.1}) and (\ref{c2}), we obtain formula (\ref{3.2}). The proof of Lemma \ref{lem3.1} is thus completed.\hfill$\square$

It is clear that using (\ref{3.1}) and applying Cauchy product formula of power series, we deduce
\begin{align}\label{c3}
\frac{{{{\ln }^k}\left( {1 + x} \right)}}{{1 - x}} = {\left( { - 1} \right)^k}k!\sum\limits_{n = 1}^\infty  {{\zeta _{n - 1}}\left( {\bar 1,{{\left\{ 1 \right\}}_{k - 1}}} \right){x^{n - 1}}} ,\quad x\in (-1,1),
\end{align}
\begin{align}\label{c4}
\frac{{{{\ln }^k}\left( {1 - x} \right)}}{{1 + x}} = {\left( { - 1} \right)^k}k!\sum\limits_{n = 1}^\infty  {{{\left( { - 1} \right)}^{n - 1}}{\zeta _{n - 1}}\left( {\bar 1,{{\left\{ 1 \right\}}_{k - 1}}} \right){x^{n - 1}}},\quad x\in(-1,1) .
\end{align}

The main results are the following theorems and corollary.
\begin{thm}\label{thm3.2} For integers $m,k\in\N_0$, then
\begin{align}\label{3.3}
\zeta \left( {\bar 1,{{\left\{ 1 \right\}}_m},\bar 1,{{\left\{ 1 \right\}}_{k}}} \right) = {\left( { - 1} \right)^{m + 1}}{\rm{L}}{{\rm{i}}_{k + 2,{{\left\{ 1 \right\}}_m}}}\left( {\frac{1}{2}} \right),
\end{align}
where the multiple polylogarithm function
\begin{align}\label{3.4}
{\rm{L}}{{\rm{i}}_{{s_1},{s_2}, \cdots {s_k}}}\left( x \right): = \sum\limits_{{n_1} >  \cdots  > {n_k} > 0} {\frac{{{x^{{n_1}}}}}{{n_1^{{s_1}} \cdots n_k^{{s_k}}}}}
\end{align}
is defined for positive integers $s_j$, and $x$ is a real number satisfying $0\leq x<1$. Of course, if $s_1>1$, then we can allow $x=1$.
\end{thm}
\pf  To prove the identity (\ref{3.3}), we consider the following multiple integral
\begin{align}\label{3.5}
{M_m}\left( k \right): = \int\limits_0^1 {\frac{1}{{1 + {t_1}}}d{t_1}}  \cdots \int\limits_0^{{t_{m - 1}}} {\frac{1}{{1 + {t_m}}}d{t_m}} \int\limits_0^{{t_m}} {\frac{{{{\ln }^{k+1}}\left( {1 - {t_{m + 1}}} \right)}}{{1 + {t_{m + 1}}}}} d{t_{m + 1}}.
\end{align}
By using formula (\ref{c4}) and the power series expansion
\[\frac{1}{{1 + x}} = \sum\limits_{n = 1}^\infty  {{{\left( { - 1} \right)}^{n - 1}}{x^{n - 1}}},\quad x\in (-1,1)\]
we can find that
\begin{align}\label{c5}
&\frac{{{{\ln }^{k + 1}}\left( {1 - {t_{m + 1}}} \right)}}{{\left( {1 + {t_1}} \right) \cdots \left( {1 + {t_m}} \right)\left( {1 + {t_{m + 1}}} \right)}} \nonumber\\
&= {\left( { - 1} \right)^{k + m}}\left( {k + 1} \right)!\sum\limits_{{n_1},{n_2}, \cdots ,{n_{m + 1}} = 1}^\infty  {{{\left( { - 1} \right)}^{{n_1} + {n_2} +  \cdots  + {n_{m + 1}}}}{\zeta _{{n_1} - 1}}\left( {\bar 1,{{\left\{ 1 \right\}}_k}} \right)t_1^{{n_{m + 1}} - 1}t_2^{{n_m} - 1} \cdots t_{m + 1}^{{n_1} - 1}} .
\end{align}
Applying the iterated-integral symbol $\int\limits_0^1 {\int\limits_0^{{t_1}} { \cdots \int\limits_0^{{t_m}} {\left(  \cdot  \right)d{t_1}d{t_2} \cdots d{t_{m + 1}}} } } $ to (\ref{c5}), we obtain
\begin{align}\label{3.6}
{M_m}\left( k \right) =& {\left( { - 1} \right)^{k + m }}(k+1)!\sum\limits_{{n_1},{n_2}, \cdots ,{n_{m + 1}} = 1}^\infty  {{{\left( { - 1} \right)}^{{n_1} +  \cdots  + {n_{m + 1}}}}\frac{{{\zeta _{{n_1} - 1}}\left( {\bar 1,{{\left\{ 1 \right\}}_{k }}} \right)}}{{{n_1}}}}\nonumber \\
&\quad \quad\quad\quad\quad\quad\quad\quad\quad\times \int\limits_0^1 {t_1^{{n_{m + 1}} - 1}d{t_1}}  \cdots \int\limits_0^{{t_{m - 1}}} {t_m^{{n_1} + {n_2} - 1}d{t_m}}\nonumber  \\
 =& {\left( { - 1} \right)^{k + m }}(k+1)!\sum\limits_{{n_1},{n_2}, \cdots ,{n_{m + 1}} = 1}^\infty  {{{\left( { - 1} \right)}^{{n_1} +  \cdots  + {n_{m + 1}}}}\frac{{{\zeta _{{n_1} - 1}}\left( {\bar 1,{{\left\{ 1 \right\}}_{k }}} \right)}}{{{n_1}\left( {{n_1} + {n_2}} \right) \cdots \left( {{n_1} +  \cdots  + {n_{m + 1}}} \right)}}}\nonumber  \\
 =& {\left( { - 1} \right)^{k + m }}(k+1)!\sum\limits_{{n_1} >  \cdots  > {n_{m + 1}} \ge 1}^\infty  {\frac{{{\zeta _{{n_{m + 1}} - 1}}\left( {\bar 1,{{\left\{ 1 \right\}}_{k}}} \right)}}{{{n_1}{n_2} \cdots {n_{m + 1}}}}{{\left( { - 1} \right)}^{{n_1}}}}\nonumber  \\
 =& {\left( { - 1} \right)^{k + m }}(k+1)!\zeta \left( {\bar 1,{{\left\{ 1 \right\}}_m},\bar 1,{{\left\{ 1 \right\}}_{k}}} \right).
\end{align}
Hence, ${M_{m}}\left( k \right) = {\left( { - 1} \right)^{k + m}}(k+1)!\zeta \left( {\bar 1,{{\left\{ 1 \right\}}_{m }},\bar 1,{{\left\{ 1 \right\}}_{k}}} \right)$.
On the other hand, applying the change of variables
\[{t_i} \mapsto 1 - {t_{m + 2 - i}},\;i = 1,2, \cdots ,m + 1.\]
to the above multiple integral ${M_m}\left( k \right)$, and noting that the fact
\[\frac{1}{{\left( {2 - {t_1}} \right) \cdots \left( {2 - {t_m}} \right)\left( {2 - {t_{m + 1}}} \right)}} = \sum\limits_{{n_1},{n_2}, \cdots ,{n_{m + 1}} = 1}^\infty  {\frac{{t_1^{{n_{m + 1}} - 1}t_2^{{n_m} - 1} \cdots t_{m + 1}^{{n_1} - 1}}}{{{2^{{n_1} + {n_2} +  \cdots  + {n_{m + 1}}}}}}} \]
we have
\begin{align}\label{3.7}
&{M_m}\left( k \right) = \int\limits_0^1 {\frac{{{{\ln }^{k+1}}\left( {{t_1}} \right)}}{{2 - {t_1}}}d{t_1}\int\limits_0^{{t_1}} {\frac{1}{{2 - {t_2}}}d{t_2}}  \cdots \int\limits_0^{{t_m}} {\frac{1}{{2 - {t_{m + 1}}}}d{t_{m + 1}}} } \nonumber \\
&{\rm{ = }}\sum\limits_{{n_1}, \cdots ,{n_{m + 1}} = 1}^\infty  {\frac{1}{{{2^{{n_1} +  \cdots  + {n_{m + 1}}}}}}\int\limits_0^1 {t_1^{{n_{m + 1}} - 1}{{\ln }^{k+1}}\left( {{t_1}} \right)d{t_1}\int\limits_0^{{t_1}} {t_2^{{n_m} - 1}d{t_2}}  \cdots \int\limits_0^{{t_m}} {t_{m + 1}^{{n_1} - 1}d{t_{m + 1}}} } }\nonumber  \\
&{\rm{ = }}\sum\limits_{{n_1}, \cdots ,{n_{m + 1}} = 1}^\infty  {\frac{1}{{{2^{{n_1} +  \cdots  + {n_{m + 1}}}}}}\cdot\frac{1}{{{n_1}\left( {{n_1} + {n_2}} \right) \cdots \left( {{n_1} +  \cdots  + {n_m}} \right)}}\int\limits_0^1 {t_1^{{n_1} +  \cdots  + {n_{m + 1}} - 1}{{\ln }^{k+1}}\left( {{t_1}} \right)d{t_1}} } \nonumber \\
& = (k+1)!{\left( { - 1} \right)^{k+1}}\sum\limits_{{n_1}, \cdots ,{n_{m + 1}} = 1}^\infty  {\frac{1}{{{2^{{n_1} +  \cdots  + {n_{m + 1}}}}}}}  \cdot \frac{1}{{{n_1}\left( {{n_1} + {n_2}} \right) \cdots \left( {{n_1} +  \cdots +{n_m}} \right){{\left( {{n_1} +  \cdots  + {n_{m + 1}}} \right)}^{k + 2}}}}\nonumber \\
&= (k+1)!{\left( { - 1} \right)^{k+1}}\sum\limits_{1 \le {n_{m + 1}} <  \cdots  < {n_1}}^\infty  {\frac{1}{{{n_{m + 1}} \cdots {n_2}n_1^{k + 2}{2^{{n_1}}}}}} \nonumber \\
& = (k+1)!{\left( { - 1} \right)^{k+1}}{\rm{L}}{{\rm{i}}_{k + 2,{{\left\{ 1 \right\}}_m}}}\left( {\frac{1}{2}} \right).
\end{align}
Therefore, combining (\ref{3.6}) and (\ref{3.7}) we may deduce the desired result.\hfill$\square$
\begin{thm}\label{thm3.3} For integers $m,k\in\N_0$, then
\begin{align}\label{3.8}
\zeta \left( {\bar 1,{{\left\{ 1 \right\}}_m},\bar 1,{{\left\{ 1 \right\}}_{k}}} \right) =& \frac{{{{\left( { - 1} \right)}^{m + k+1}}}}{{ {k} !}}\sum\limits_{j = 0}^{k} {{{\left( { - 1} \right)}^j}{{\left( {\ln 2} \right)}^{k- j}}j!\left( {\begin{array}{*{20}{c}}
   {k}  \\
   j  \\
\end{array}} \right)}\nonumber  \\
&\quad\quad\quad\quad\times \left\{ {\zeta \left( {m + 2,{{\left\{ 1 \right\}}_j}} \right) - \sum\limits_{l = 0}^{m + 1} {\frac{{{{\left( {\ln 2} \right)}^{m + 1 - l}}}}{{\left( {m + 1 - l} \right)!}}{\rm{L}}{{\rm{i}}_{l + 1,{{\left\{ 1 \right\}}_j}}}\left( {\frac{1}{2}} \right)} } \right\}.
\end{align}
\end{thm}
\pf By using (\ref{3.1}), we can find that
\begin{align}\label{3.9}
\int\limits_0^1 {\frac{{{{\left( {\ln x} \right)}^k}{{\ln }^{m + 1}}\left( {1 - \frac{x}{2}} \right)}}{x}} dx &= {\left( { - 1} \right)^{m + 1}}\left( {m + 1} \right)!\sum\limits_{n = 1}^\infty  {\frac{{{\zeta _{n - 1}}\left( {{{\left\{ 1 \right\}}_m}} \right)}}{{n{2^n}}}\int\limits_0^1 {{x^{n - 1}}{{\ln }^k}(x)dx} }\nonumber  \\
&= {\left( { - 1} \right)^{m + k + 1}}\left( {m + 1} \right)!k!\sum\limits_{n = 1}^\infty  {\frac{{{\zeta _{n - 1}}\left( {{{\left\{ 1 \right\}}_m}} \right)}}{{{n^{k + 2}}{2^n}}}} \nonumber \\
& = {\left( { - 1} \right)^{m + k + 1}}\left( {m + 1} \right)!k!{\rm{L}}{{\rm{i}}_{k + 2,{{\left\{ 1 \right\}}_m}}}\left( {\frac{1}{2}} \right).
\end{align}
Then it is readily seen that
\begin{align}\label{3.10}
{\rm{L}}{{\rm{i}}_{k + 2,{{\left\{ 1 \right\}}_m}}}\left( {\frac{1}{2}} \right) = \frac{{{{\left( { - 1} \right)}^{m + k+1}}}}{{\left( {m + 1} \right)!{k} !}}\int\limits_0^1 {\frac{{{{\left( {\ln x} \right)}^{k}}{{\ln }^{m + 1}}\left( {1 - \frac{x}{2}} \right)}}{x}} dx,
\end{align}
\begin{align}\label{3.11}
\zeta \left( {\bar 1,{{\left\{ 1 \right\}}_m},\bar 1,{{\left\{ 1 \right\}}_{k}}} \right) = \frac{{{{\left( { - 1} \right)}^{k}}}}{{\left( {m + 1} \right)! {k}!}}\int\limits_0^1 {\frac{{{{\left( {\ln x} \right)}^{k}}{{\ln }^{m + 1}}\left( {1 - \frac{x}{2}} \right)}}{x}} dx.
\end{align}
On the other hand, changing $x=2(1-u)$ in above integral, we conclude that
\begin{align}\label{3.12}
&\int\limits_0^1 {\frac{{{{\left( {\ln x} \right)}^{k}}{{\ln }^{m + 1}}\left( {1 - \frac{x}{2}} \right)}}{x}} dx\nonumber \\&= \int\limits_{\frac{1}{2}}^1 {\frac{{{{\left( {\ln 2 + \ln \left( {1 - u} \right)} \right)}^{k}}{{\ln }^{m + 1}}(u)}}{{1 - u}}} du\nonumber\\
& = \sum\limits_{j = 0}^{k} {{{\left( {\ln 2} \right)}^{k- j}}\left( {\begin{array}{*{20}{c}}
   {k}  \\
   j  \\
\end{array}} \right)} \int\limits_{\frac{1}{2}}^1 {\frac{{{{\ln }^j}\left( {1 - u} \right){{\ln }^{m + 1}}(u)}}{{1 - u}}} du\nonumber\\
& = \sum\limits_{j = 0}^{k} {{{\left( { - 1} \right)}^j}{{\left( {\ln 2} \right)}^{k - j}}j!\left( {\begin{array}{*{20}{c}}
   {k}  \\
   j  \\
\end{array}} \right)\sum\limits_{n = 1}^\infty  {{\zeta _{n - 1}}\left( {{{\left\{ 1 \right\}}_j}} \right)\int\limits_{\frac{1}{2}}^1 {{u^{n - 1}}{{\ln }^{m + 1}}(u)} du} }\nonumber \\
& = \left( {m + 1} \right)!{\left( { - 1} \right)^{m + 1}}\sum\limits_{j = 0}^{k} {{{\left( { - 1} \right)}^j}{{\left( {\ln 2} \right)}^{k - j}}j!\left( {\begin{array}{*{20}{c}}
   {k}  \\
   j  \\
\end{array}} \right)}\nonumber \\
&\quad\quad\quad\quad\quad\quad\quad\quad\quad\times \left\{ {\zeta \left( {m + 2,{{\left\{ 1 \right\}}_j}} \right) - \sum\limits_{l = 0}^{m + 1} {\frac{{{{\left( {\ln 2} \right)}^{m + 1 - l}}}}{{\left( {m + 1 - l} \right)!}}{\rm{L}}{{\rm{i}}_{l + 1,{{\left\{ 1 \right\}}_j}}}\left( {\frac{1}{2}} \right)} } \right\}.
\end{align}
Thus, substituting (\ref{3.12}) into (\ref{3.11}), we obtain the desired result. The proof of Theorem \ref{thm3.3} is finished. \hfill$\square$\\
From Theorem \ref{thm3.2} and Theorem \ref{thm3.3}, we immediately derive the following special case of the multiple zeta values.
\begin{cor} (Conjectured in \cite{BBBL1997})\label{cor3.4} For integer $m\in\N_0$, then the following identities hold:
\begin{align}\label{3.13}
&{\rm{L}}{{\rm{i}}_{2,{{\left\{ 1 \right\}}_m}}}\left( {\frac{1}{2}} \right) = \zeta \left( {m + 2} \right) - \sum\limits_{l = 0}^{m + 1} {\frac{{{{\left( {\ln 2} \right)}^{m + 1 - l}}}}{{\left( {m + 1 - l} \right)!}}{\rm{L}}{{\rm{i}}_{l + 1}}\left( {\frac{1}{2}} \right)} ,
\end{align}
\begin{align}\label{3.14}
\zeta \left( {\bar 1,{{\left\{ 1 \right\}}_m},\bar 1} \right) = {\left( { - 1} \right)^{m + 1}}\left\{ {\zeta \left( {m + 2} \right) - \sum\limits_{l = 0}^{m + 1} {\frac{{{{\left( {\ln 2} \right)}^{m + 1 - l}}}}{{\left( {m + 1 - l} \right)!}}{\rm{L}}{{\rm{i}}_{l + 1}}\left( {\frac{1}{2}} \right)} } \right\}.
\end{align}
\end{cor}
\pf Setting $k=0$ in (\ref{3.3}) and (\ref{3.8}) we may easily deduce the results.\hfill$\square$\\
Note that Corollary \ref{cor3.4} is an immediate corollary of Zlobin's Theorem 9 (see \cite{SAZ2012}).
\begin{thm}\label{thm3.5} For integers $m,k\in \N_0$, we have
\begin{align}\label{3.15}
\zeta \left( {\bar 1,{{\left\{ 1 \right\}}_{m}},\bar 1,\bar 1,{{\left\{ 1 \right\}}_{k}}} \right) =& \frac{{{{\left( { - 1} \right)}^{k}}}}{{(m+1)!(k+1)!}}\left\{ {(k+1){{\left( {\ln 2} \right)}^{m+1}}I\left( {k} \right) - \left( {m + k+2} \right)I\left( {m + k+ 1} \right)} \right\}\nonumber\\
& - \sum\limits_{i = 1}^{m} {\frac{{{{\left( {\ln 2} \right)}^i}}}{{i!}}} \zeta \left( {\bar 1,{{\left\{ 1 \right\}}_{m - i}},\bar 1,\bar 1,{{\left\{ 1 \right\}}_{k}}} \right),
\end{align}
where $I\left( k \right)$ denotes the integral on the left hand side of (\ref{2.4}), namely
\[I\left( k \right): = \int\limits_0^1 {\frac{{{{\ln }^k}\left( {1 + x} \right)\ln \left( {1 - x} \right)}}{{1 + x}}dx}\quad (k\in \N),\]
with \[I\left( 0 \right) = \frac{{{{\ln }^2}\left( 2 \right) - \zeta \left( 2 \right)}}{2}.\]
\end{thm}
\pf By a similar argument as in the proof of Theorem \ref{thm3.2}, we consider the following multiple integral
\begin{align}\label{3.16}
{J_m}\left( k \right): = \int\limits_0^1 {\frac{1}{{1 + {t_1}}}d{t_1} \cdots } \int\limits_0^{{t_{m - 2}}} {\frac{1}{{1 + {t_{m - 1}}}}d{t_{m - 1}}} \int\limits_0^{{t_{m - 1}}} {\frac{1}{{1 + {t_m}}}d{t_m}} \int\limits_0^{{t_m}} {\frac{{{{\ln }^k}\left( {1 + {t_{m + 1}}} \right)}}{{1 - {t_{m + 1}}}}d{t_{m + 1}}}.
\end{align}
By using (\ref{c3}), we deduce that
\begin{align}\label{3.17}
\int\limits_0^x {\frac{1}{{1 + t}}dt\int\limits_0^t {\frac{{{{\ln }^k}\left( {1 + u} \right)}}{{1 - u}}du} }  = {\left( { - 1} \right)^{k - 1}}k!\sum\limits_{n = 1}^\infty  {\frac{{{\zeta _{n - 1}}\left( {\bar 1,\bar 1,{{\left\{ 1 \right\}}_{k - 1}}} \right)}}{n}{{\left( { - 1} \right)}^n}{x^n}} .
\end{align}
Hence, combining (\ref{3.16}) and (\ref{3.17}), it is easily shown that
\begin{align}\label{3.18}
{J_m}\left( k \right) = {\left( { - 1} \right)^{m + k}}k! \zeta \left( {\bar 1,{{\left\{ 1 \right\}}_{m - 1}},\bar 1,\bar 1,{{\left\{ 1 \right\}}_{k - 1}}} \right).
\end{align}
On the other hand, using integration by parts, we can find that
\begin{align}\label{c6}
 &{J_m}\left( k \right) = \int\limits_0^1 {\frac{1}{{1 + {t_1}}}\left( {\int\limits_{0 < {t_{m + 1}} <  \cdots  < {t_2} < {t_1}} {\frac{{{{\ln }^k}\left( {1 + {t_{m + 1}}} \right)}}{{\left( {1 + {t_2}} \right) \cdots \left( {1 + {t_m}} \right)\left( {1 - {t_{m + 1}}} \right)}}d{t_2} \cdots d{t_{m + 1}}} } \right)d{t_1}}\nonumber  \\
  =& \left. {\left( {\ln \left( {1 + {t_1}} \right)\int\limits_{0 < {t_{m + 1}} <  \cdots  < {t_2} < {t_1}} {\frac{{{{\ln }^k}\left( {1 + {t_{m + 1}}} \right)}}{{\left( {1 + {t_2}} \right) \cdots \left( {1 + {t_m}} \right)\left( {1 - {t_{m + 1}}} \right)}}d{t_2} \cdots d{t_{m + 1}}} } \right)} \right|_{{t_1} = 0}^{{t_1} = 1} \nonumber \\
 & - \int\limits_0^1 {\frac{{\ln \left( {1 + {t_1}} \right)}}{{1 + {t_1}}}} \int\limits_{0 < {t_{m + 1}} <  \cdots  < {t_3} < {t_1}} {\frac{{{{\ln }^k}\left( {1 + {t_{m + 1}}} \right)}}{{\left( {1 + {t_3}} \right) \cdots \left( {1 + {t_m}} \right)\left( {1 - {t_{m + 1}}} \right)}}d{t_1}d{t_3} \cdots d{t_{m + 1}}} \nonumber  \\
  = &\left( {\ln 2} \right){J_{m - 1}}\left( k \right) - \int\limits_0^1 {\frac{{\ln \left( {1 + {t_1}} \right)}}{{1 + {t_1}}}d{t_1}\int\limits_0^{{t_1}} {\frac{1}{{1 + {t_2}}}d{t_2}}  \cdots \int\limits_0^{{t_{m - 2}}} {\frac{1}{{1 + {t_{m - 1}}}}d{t_{m - 1}}\int\limits_0^{{t_{m - 1}}} {\frac{{{{\ln }^k}\left( {1 + {t_m}} \right)}}{{1 - {t_m}}}d{t_m}} } } \nonumber\\
 =&\cdots\nonumber\\
 =& \sum\limits_{i = 1}^{m - 1} {{{\left( { - 1} \right)}^{i - 1}}\frac{{{{\left( {\ln 2} \right)}^i}}}{{i!}}{J_{m - i}}\left( k \right)}  + \frac{{{{\left( { - 1} \right)}^{m - 1}}}}{{\left( {m - 1} \right)!}}\int\limits_0^1 {\frac{{{{\ln }^{m - 1}}\left( {1 + {t_1}} \right)}}{{1 + {t_1}}}d{t_1}\int\limits_0^{{t_1}} {\frac{{{{\ln }^k}\left( {1 + {t_2}} \right)}}{{1 - {t_2}}}d{t_2}} } .
\end{align}
Therefore, substituting (\ref{3.18}) into (\ref{c6}) yields
\begin{align}\label{3.19}
{J_m}\left( k \right) = &{\left( { - 1} \right)^{m + k - 1}}k!\sum\limits_{i = 1}^{m - 1} {\frac{{{{\left( {\ln 2} \right)}^i}}}{{i!}}} \zeta \left( {\bar 1,{{\left\{ 1 \right\}}_{m - i - 1}},\bar 1,\bar 1,{{\left\{ 1 \right\}}_{k - 1}}} \right)\nonumber\\
& + \frac{{{{\left( { - 1} \right)}^{m - 1}}}}{{\left( {m - 1} \right)!}}\int\limits_0^1 {\frac{{{{\ln }^{m - 1}}\left( {1 + {t_1}} \right)}}{{1 + {t_1}}}d{t_1}\int\limits_0^{{t_1}} {\frac{{{{\ln }^k}\left( {1 + {t_2}} \right)}}{{1 - {t_2}}}d{t_2}} } .
\end{align}
Moreover, we note that the integral on the right-hand side of (\ref{3.19}) can be written as
\begin{align}\label{3.20}
&\int\limits_0^1 {\frac{{{{\ln }^{m - 1}}\left( {1 + {t_1}} \right)}}{{1 + {t_1}}}d{t_1}\int\limits_0^{{t_1}} {\frac{{{{\ln }^k}\left( {1 + {t_2}} \right)}}{{1 - {t_2}}}d{t_2}} }\nonumber \\
& = \mathop {\lim }\limits_{x \to 1} \left\{ {\int\limits_0^x {\frac{{{{\ln }^{m - 1}}\left( {1 + {t_1}} \right)}}{{1 + {t_1}}}d{t_1}\int\limits_0^{{t_1}} {\frac{{{{\ln }^k}\left( {1 + {t_2}} \right)}}{{1 - {t_2}}}d{t_2}} } } \right\}\nonumber\\
& = \frac{1}{m}\mathop {\lim }\limits_{x \to 1} \left\{ {\int\limits_0^x {\frac{{{{\ln }^m}\left( {1 + x} \right){{\ln }^k}\left( {1 + t} \right) - {{\ln }^{m + k}}\left( {1 + t} \right)}}{{1 - t}}dt} } \right\}\nonumber\\
& = \frac{1}{m}\left\{ {k{{\left( {\ln 2} \right)}^m}\int\limits_0^1 {\frac{{{{\ln }^{k - 1}}\left( {1 + t} \right)\ln \left( {1 - t} \right)}}{{1 + t}}dt}  - \left( {m + k} \right)\int\limits_0^1 {\frac{{{{\ln }^{m + k - 1}}\left( {1 + t} \right)\ln \left( {1 - t} \right)}}{{1 + t}}dt} } \right\}.
\end{align}
Therefore, replacing $k-1$ by $k$ and $m-1$ by $m$, the relations (\ref{3.18}), (\ref{3.19}) and (\ref{3.20}) yield the desired result.  The proof of Theorem 3.5 is completed.\hfill$\square$

From Lemma \ref{lem2.2} and Theorem \ref{thm3.5}, we have the conclusion: if $m,k\in \N_0$, then the alternating multiple zeta values $\zeta \left( {\bar 1,{{\left\{ 1 \right\}}_{m}},\bar 1,\bar 1,{{\left\{ 1 \right\}}_{k}}} \right)$
can be expressed as a rational linear combination of zeta values, polylogarithms and $\ln(2)$.
We now close this section with several examples.
\begin{exa} By (\ref{3.14}) and (\ref{3.15}), we have
\begin{align*}
&\zeta \left( {\bar 1,1,\bar 1} \right) = \frac{1}{8}\zeta \left( 3 \right) - \frac{1}{6}{\ln ^3}(2),\\
&\zeta \left( {\bar 1,\bar1,\bar 1} \right) = -\frac{1}{4}\zeta \left( 3 \right)+\frac{1}{2}\zeta \left( 2 \right)\ln(2) - \frac{1}{6}{\ln ^3}(2),\\
&\zeta \left( {\bar 1,1,1,\bar 1} \right) = {\rm{L}}{{\rm{i}}_4}\left( {\frac{1}{2}} \right) + \frac{1}{{12}}{\ln ^4}(2) + \frac{7}{8}\zeta \left( 3 \right)\ln(2) - \frac{1}{2}\zeta \left( 2 \right){\ln ^2}(2) - \zeta \left( 4 \right),\\
&\zeta \left( {\bar 1,\bar 1,\bar 1,1} \right) = 3{\rm{L}}{{\rm{i}}_4}\left( {\frac{1}{2}} \right) + \frac{1}{6}{\ln ^4}(2) + \frac{{23}}{8}\zeta \left( 3 \right)\ln (2) - \zeta \left( 2 \right){\ln ^2}(2) - 3\zeta \left( 4 \right),\\
&\zeta \left( {\bar 1,1,\bar 1,\bar 1} \right) =  - 3{\rm{L}}{{\rm{i}}_4}\left( {\frac{1}{2}} \right) - \frac{1}{{12}}{\ln ^4}(2)- \frac{{11}}{4}\zeta \left( 3 \right)\ln (2) - \frac{3}{4}\zeta \left( 2 \right){\ln ^2}(2) + 3\zeta \left( 4 \right).
\end{align*}
\end{exa}
\section{Some evaluation of restricted sum formulas involving multiple zeta values}
In \cite{X2017}, we considered the following restricted sum formulas involving multiple zeta values
\[{\left( { - 1} \right)^{m+1}}\sum\limits_{i = 0}^{m} {\frac{{{{\left( {\ln 2} \right)}^i}}}{{i!}}\zeta \left( {\bar 1,{{\left\{ 1 \right\}}_{m - i}},3,{{\left\{ 1 \right\}}_{k}}} \right)}  + {\left( { - 1} \right)^{k+1}}\sum\limits_{i = 0}^{k} {\frac{{{{\left( {\ln 2} \right)}^i}}}{{i!}}\zeta \left( {\bar 1,{{\left\{ 1 \right\}}_{k - i}},3,{{\left\{ 1 \right\}}_{m }}} \right)} \]
and gave explicit reductions to multiple zeta values of depth less than ${\rm max}\{k+3,m+3\}$, where $m,k\in \N_0$. Moreover, we also proved that the restricted sum
\[\sum\limits_{i = 0}^{m} {\frac{{{{\left( {\ln 2} \right)}^i}}}{{i!}}\zeta \left( {\bar 1,{{\left\{ 1 \right\}}_{m - i}},2,{{\left\{ 1 \right\}}_{k}}} \right)} \]
can be expressed by the zeta values and polylogarithms, which implies that for any $m,k\in \N_0$, the multiple zeta values $\zeta \left( {\bar 1,{{\left\{ 1 \right\}}_{m}},2,{{\left\{ 1 \right\}}_{k}}} \right)$ can be represented as a polynomial of zeta values and polylogarithms with rational coefficients. In particular, one can find explicit formula for weight 4
\begin{align*}
\zeta \left( {\bar 1,1,2} \right) = 3{\rm{L}}{{\rm{i}}_4}\left( {\frac{1}{2}} \right) + \frac{1}{8}{\ln ^4}(2) + \frac{{23}}{8}\zeta \left( 3 \right)\ln (2) - \zeta \left( 2 \right){\ln ^2}(2) - 3\zeta \left( 4 \right).
\end{align*}
In this section, we will consider the general restricted sum
\[{\left( { - 1} \right)^{m+1}}\sum\limits_{i = 0}^{m} {\frac{{{{\left( {\ln 2} \right)}^i}}}{{i!}}\zeta \left( {\bar 1,{{\left\{ 1 \right\}}_{m - i}},p + 3,{{\left\{ 1 \right\}}_{k}}} \right)}  + {\left( { - 1} \right)^{p + k+1}}\sum\limits_{i = 0}^{k} {\frac{{{{\left( {\ln 2} \right)}^i}}}{{i!}}\zeta \left( {\bar 1,{{\left\{ 1 \right\}}_{k - i}},p + 3,{{\left\{ 1 \right\}}_{m }}} \right)} \]
where $m,k,p\in\N_0$.
Now we are ready to state and prove our main results. Note that our proof of Theorem \ref{thm4.1} is based on Lemma \ref{lem3.1}.
\begin{thm}\label{thm4.1} For integers $m,k,p\in\N_0$, we have
\begin{align}\label{4.1}
&{\left( { - 1} \right)^{m+1}}\sum\limits_{i = 0}^{m} {\frac{{{{\left( {\ln 2} \right)}^i}}}{{i!}}\zeta \left( {\bar 1,{{\left\{ 1 \right\}}_{m - i}},p + 3,{{\left\{ 1 \right\}}_{k}}} \right)}\nonumber \\
&\quad + {\left( { - 1} \right)^{p + k+1}}\sum\limits_{i = 0}^{k} {\frac{{{{\left( {\ln 2} \right)}^i}}}{{i!}}\zeta \left( {\bar 1,{{\left\{ 1 \right\}}_{k - i}},p + 3,{{\left\{ 1 \right\}}_{m}}} \right)} \nonumber\\
&= \frac{{{{\left( { - 1} \right)}^{m }}}}{{(m+1)!}}{\left( {\ln 2} \right)^{m+1}}\zeta \left( {{\overline {p + 3}},{{\left\{ 1 \right\}}_{k}}} \right)\nonumber\\
&\quad + \frac{{{{\left( { - 1} \right)}^{p + k }}}}{{(k+1)!}}{\left( {\ln 2} \right)^{k+1}}\zeta \left( {{\overline {p + 3}},{{\left\{ 1 \right\}}_{m}}} \right)\nonumber\\
&\quad + \sum\limits_{i = 0}^p {{{\left( { - 1} \right)}^i}\zeta \left( {{\overline {2 + i}},{{\left\{ 1 \right\}}_{m}}} \right)} \zeta \left( {{\overline {p + 2 - i}},{{\left\{ 1 \right\}}_{k }}} \right).
\end{align}
\end{thm}
\pf Similarly as in the proof of Theorem \ref{thm3.5}, we consider the multiple integral
\[
{R_{m,k}}\left( p \right): = \int\limits_{0 < {t_{m + p + 2}} <  \cdots  < {t_1} < 1} {\frac{{{{\ln }^k}\left( {1 + {t_{m + p + 2}}} \right)}}{{\left( {1 + {t_1}} \right) \cdots \left( {1 + {t_m}} \right){t_{m + 1}} \cdots {t_{m + p + 1}}{t_{m + p + 2}}}}} d{t_1} \cdots d{t_{m + p + 2}}.
\]
Then with the help of formula (\ref{3.1}), we may easily deduce that
\begin{align}\label{4.2}
{R_{m,k}}\left( p \right) = {\left( { - 1} \right)^{m + k}}k!\zeta \left( {\bar 1,{{\left\{ 1 \right\}}_{m - 1}},p + 3,{{\left\{ 1 \right\}}_{k - 1}}} \right).
\end{align}
On the other hand, by a similar argument as in the proof of formula (\ref{c6}) (using integration by parts), it is easily shown that
\begin{align}\label{4.3}
&{\left( { - 1} \right)^{m + k}}k!\zeta \left( {\bar 1,{{\left\{ 1 \right\}}_{m - 1}},p + 3,{{\left\{ 1 \right\}}_{k - 1}}} \right)\nonumber\\
& = {\left( { - 1} \right)^{m + k + 1}}k!\sum\limits_{i = 1}^{m - 1} {\frac{{{{\left( {\ln 2} \right)}^i}}}{{i!}}\zeta \left( {\bar 1,{{\left\{ 1 \right\}}_{m - 1 - i}},p + 3,{{\left\{ 1 \right\}}_{k - 1}}} \right)} \nonumber\\
&\quad + \frac{{{{\left( { - 1} \right)}^{m - 1}}}}{{\left( {m - 1} \right)!}}\int\limits_0^1 {\frac{{{{\ln }^{m - 1}}\left( {1 + {t_1}} \right)}}{{1 + {t_1}}}d{t_1}} \int\limits_0^{{t_1}} {\frac{1}{{{t_2}}}d{t_2} \cdots \int\limits_0^{{t_{p{\rm{ + }}1}}} {\frac{1}{{{t_{p + 2}}}}{t_{p + 2}}} } \int\limits_0^{{t_{p + 2}}} {\frac{{{{\ln }^k}\left( {1 + {t_{p + 3}}} \right)}}{{{t_{p + 3}}}}d{t_{p + 3}}}\nonumber \\
& = {\left( { - 1} \right)^{m + k + 1}}k!\sum\limits_{i = 1}^{m - 1} {\frac{{{{\left( {\ln 2} \right)}^i}}}{{i!}}\zeta \left( {\bar 1,{{\left\{ 1 \right\}}_{m - 1 - i}},p + 3,{{\left\{ 1 \right\}}_{k - 1}}} \right)}\nonumber \\
&\quad + \frac{{{{\left( { - 1} \right)}^{m +k- 1}}}}{{m!}}{k!\left( {\ln 2} \right)^m}\zeta \left( {{\overline {p + 3}},{{\left\{ 1 \right\}}_{k - 1}}} \right)\nonumber\\
&\quad + \frac{{{{\left( { - 1} \right)}^m}}}{{m!}}\int\limits_0^1 {\frac{{{{\ln }^m}\left( {1 + {t_1}} \right)}}{{{t_1}}}d{t_1}} \int\limits_0^{{t_1}} {\frac{1}{{{t_2}}}d{t_2} \cdots \int\limits_0^{{t_p}} {\frac{1}{{{t_{p + 1}}}}{t_{p + 1}}} } \int\limits_0^{{t_{p + 1}}} {\frac{{{{\ln }^k}\left( {1 + {t_{p + 2}}} \right)}}{{{t_{p + 2}}}}d{t_{p + 2}}} .
\end{align}
Hence, by a direct calculation, we can get the following identity
\begin{align}\label{4.4}
&\int\limits_0^1 {\frac{{{{\ln }^m}\left( {1 + {t_1}} \right)}}{{{t_1}}}d{t_1}} \int\limits_0^{{t_1}} {\frac{1}{{{t_2}}}d{t_2} \cdots \int\limits_0^{{t_p}} {\frac{1}{{{t_{p + 1}}}}{t_{p + 1}}} } \int\limits_0^{{t_{p + 1}}} {\frac{{{{\ln }^k}\left( {1 + {t_{p + 2}}} \right)}}{{{t_{p + 2}}}}d{t_{p + 2}}}\nonumber \\
& = {\left( { - 1} \right)^k}m!k!\sum\limits_{i = 0}^{m - 1} {\frac{{{{\left( {\ln 2} \right)}^i}}}{{i!}}\zeta \left( {\bar 1,{{\left\{ 1 \right\}}_{m - 1 - i}},p + 3,{{\left\{ 1 \right\}}_{k - 1}}} \right)}\nonumber \\
&\quad + {\left( { - 1} \right)^k}k!{\left( {\ln 2} \right)^m}\zeta \left( {{\overline {p + 3}},{{\left\{ 1 \right\}}_{k - 1}}} \right).
\end{align}
Moreover, applying the same argument as above, we deduce that
\begin{align}\label{4.5}
&\int\limits_0^1 {\frac{{{{\ln }^m}\left( {1 + {t_1}} \right)}}{{{t_1}}}d{t_1}} \int\limits_0^{{t_1}} {\frac{1}{{{t_2}}}d{t_2} \cdots \int\limits_0^{{t_p}} {\frac{1}{{{t_{p + 1}}}}{t_{p + 1}}} } \int\limits_0^{{t_{p + 1}}} {\frac{{{{\ln }^k}\left( {1 + {t_{p + 2}}} \right)}}{{{t_{p + 2}}}}d{t_{p + 2}}}\nonumber \\
& + {\left( { - 1} \right)^p}\int\limits_0^1 {\frac{{{{\ln }^k}\left( {1 + {t_1}} \right)}}{{{t_1}}}d{t_1}} \int\limits_0^{{t_1}} {\frac{1}{{{t_2}}}d{t_2} \cdots \int\limits_0^{{t_p}} {\frac{1}{{{t_{p + 1}}}}{t_{p + 1}}} } \int\limits_0^{{t_{p + 1}}} {\frac{{{{\ln }^m}\left( {1 + {t_{p + 2}}} \right)}}{{{t_{p + 2}}}}d{t_{p + 2}}}\nonumber \\
&= {\left( { - 1} \right)^{m + k}}k!m!\sum\limits_{i = 0}^p {{{\left( { - 1} \right)}^i}\zeta \left( {{\overline {2 + i}},{{\left\{ 1 \right\}}_{m - 1}}} \right)} \zeta \left( {{\overline {p + 2 - i}},{{\left\{ 1 \right\}}_{k - 1}}} \right).
\end{align}
Thus, combining the formulas (\ref{4.4}) and (\ref{4.5}), then replacing $k-1$ by $k$ and $m-1$ by $m$, we obtain the desired result. This completes the proof of Theorem \ref{thm4.1}.\hfill$\square$\\
Taking $p=0$ in Theorem \ref{thm4.1}, we get the following corollary which was first proved in \cite{X2017}.
\begin{cor} For integers $m,k\in \N_0$, we have
\begin{align}\label{4.6}
&{\left( { - 1} \right)^{m+1}}\sum\limits_{i = 0}^{m } {\frac{{{{\left( {\ln 2} \right)}^i}}}{{i!}}\zeta \left( {\bar 1,{{\left\{ 1 \right\}}_{m - i}},3,{{\left\{ 1 \right\}}_{k}}} \right)}\nonumber \\
&+ {\left( { - 1} \right)^{k+1}}\sum\limits_{i = 0}^{k} {\frac{{{{\left( {\ln 2} \right)}^i}}}{{i!}}\zeta \left( {\bar 1,{{\left\{ 1 \right\}}_{k - i}},3,{{\left\{ 1 \right\}}_{m}}} \right)}\nonumber \\
& = \frac{{{{\left( { - 1} \right)}^{m}}}}{{(m+1)!}}{\left( {\ln 2} \right)^{m+1}}\zeta \left( {\bar 3,{{\left\{ 1 \right\}}_{k}}} \right) + \frac{{{{\left( { - 1} \right)}^{k  }}}}{{(k+1)!}}{\left( {\ln 2} \right)^{k+1}}\zeta \left( {\bar 3,{{\left\{ 1 \right\}}_{m}}} \right)\nonumber\\
&\quad + \zeta \left( {\bar 2,{{\left\{ 1 \right\}}_{m}}} \right)\zeta \left( {\bar 2,{{\left\{ 1 \right\}}_{k}}} \right).
\end{align}
\end{cor}
\begin{thm}\label{thm4.3} For integers $m,k,p\in \N_0$, then the following identity holds:
\begin{align}\label{4.7}
&{\left( { - 1} \right)^{k + p}}(m+1)!(k+1)!\sum\limits_{i = 0}^{m} {\frac{{{{\left( {\ln 2} \right)}^i}}}{{i!}}\zeta \left( {\bar 1,{{\left\{ 1 \right\}}_{m - i}},2,{{\left\{ 1 \right\}}_{p}},2,{{\left\{ 1 \right\}}_{k}}} \right)}\nonumber  \\
&+ {\left( { - 1} \right)^{m+1}}(k+1)!(m+1)!\sum\limits_{i = 0}^{k} {\frac{{{{\left( {\ln 2} \right)}^i}}}{{i!}}\zeta \left( {\bar 1,{{\left\{ 1 \right\}}_{k - i}},2,{{\left\{ 1 \right\}}_{p}},2,{{\left\{ 1 \right\}}_{m}}} \right)}\nonumber  \\
&+ {\left( { - 1} \right)^{k + p}}(k+1)!{\left( {\ln 2} \right)^{m+1}}\zeta \left( {\bar 2,{{\left\{ 1 \right\}}_{p}},2,{{\left\{ 1 \right\}}_{k}}} \right)\nonumber \\
& + {\left( { - 1} \right)^{m+1}}(m+1)!{\left( {\ln 2} \right)^{k+1}}\zeta \left( {\bar 2,{{\left\{ 1 \right\}}_{p}},2,{{\left\{ 1 \right\}}_{m}}} \right)\nonumber \\
& = {\left( { - 1} \right)^{m + k + p+1}}{(m+1)}!{(k+1)}!\zeta \left( {\bar 2,{{\left\{ 1 \right\}}_{m}}} \right)\zeta \left( {\bar 1,{{\left\{ 1 \right\}}_{p}},2,{{\left\{ 1 \right\}}_{k}}} \right)\nonumber \\
&\quad + {\left( { - 1} \right)^{m + k}}(k+1)!(m+1)!\zeta \left( {\bar 2,{{\left\{ 1 \right\}}_{k}}} \right)\zeta \left( {\bar 1,{{\left\{ 1 \right\}}_{p}},2,{{\left\{ 1 \right\}}_{m }}} \right)\nonumber \\
&\quad + {\left( { - 1} \right)^{m + k + p+1}}(k+1)!(m+1)!\sum\limits_{i = 1}^{p} {{{\left( { - 1} \right)}^i}\zeta \left( {\bar 1,{{\left\{ 1 \right\}}_{i - 1}},2,{{\left\{ 1 \right\}}_{m }}} \right)\zeta \left( {\bar 1,{{\left\{ 1 \right\}}_{p - i}},2,{{\left\{ 1 \right\}}_{k}}} \right)} .
\end{align}
\end{thm}
\pf The proof of Theorem \ref{thm4.3} is similar to the proof of Theorem \ref{thm4.1}. We consider the iterated integral
\[{Q_{m,k}}\left( p \right): = \int\limits_{0 < {t_{m + p + 2}} <  \cdots  < {t_1} < 1} {\frac{{{{\ln }^k}\left( {1 + {t_{m + p + 2}}} \right)d{t_1} \cdots d{t_{m + p + 2}}}}{{\left( {1 + {t_1}} \right) \cdots \left( {1 + {t_m}} \right){t_{m + 1}}\left( {1 + {t_{m + 2}}} \right) \cdots \left( {1 + {t_{m + p + 1}}} \right){t_{m + p + 2}}}}} .\]
By a similar argument as in the proof of formula (\ref{4.1}), we can obtain the following identities
\begin{align}\label{4.8}
&{Q_{m,k}}\left( p \right) = {\left( { - 1} \right)^{k + m + p}}k!\zeta \left( {\bar 1,{{\left\{ 1 \right\}}_{m - 1}},2,{{\left\{ 1 \right\}}_{p - 1}},2,{{\left\{ 1 \right\}}_{k - 1}}} \right),
\end{align}
\begin{align}\label{4.9}
&\int\limits_0^1 {\frac{{{{\ln }^m}\left( {1 + {t_1}} \right)}}{{{t_1}}}d{t_1}} \int\limits_0^{{t_1}} {\frac{1}{{1 + {t_2}}}d{t_2} \cdots \int\limits_0^{{t_p}} {\frac{1}{{1 + {t_{p + 1}}}}{t_{p + 1}}} } \int\limits_0^{{t_{p + 1}}} {\frac{{{{\ln }^k}\left( {1 + {t_{p + 2}}} \right)}}{{{t_{p + 2}}}}d{t_{p + 2}}}\nonumber  \\
& = {\left( { - 1} \right)^{k + p}}m!k!\sum\limits_{i = 0}^{m - 1} {\frac{{{{\left( {\ln 2} \right)}^i}}}{{i!}}\zeta \left( {\bar 1,{{\left\{ 1 \right\}}_{m - i - 1}},2,{{\left\{ 1 \right\}}_{p - 1}},2,{{\left\{ 1 \right\}}_{k - 1}}} \right)}\nonumber  \\
&\quad+ {\left( { - 1} \right)^{p + k}}k!{\left( {\ln 2} \right)^m}\zeta \left( {\bar 2,{{\left\{ 1 \right\}}_{p - 1}},2,{{\left\{ 1 \right\}}_{k - 1}}} \right),
\end{align}
\begin{align}\label{4.10}
&\int\limits_0^1 {\frac{{{{\ln }^m}\left( {1 + {t_1}} \right)}}{{{t_1}}}d{t_1}} \int\limits_0^{{t_1}} {\frac{1}{{1 + {t_2}}}d{t_2} \cdots \int\limits_0^{{t_p}} {\frac{1}{{1 + {t_{p + 1}}}}{t_{p + 1}}} } \int\limits_0^{{t_{p + 1}}} {\frac{{{{\ln }^k}\left( {1 + {t_{p + 2}}} \right)}}{{{t_{p + 2}}}}d{t_{p + 2}}}\nonumber  \\
&+ {\left( { - 1} \right)^p}\int\limits_0^1 {\frac{{{{\ln }^k}\left( {1 + {t_1}} \right)}}{{{t_1}}}d{t_1}} \int\limits_0^{{t_1}} {\frac{1}{{1 + {t_2}}}d{t_2} \cdots \int\limits_0^{{t_p}} {\frac{1}{{1 + {t_{p + 1}}}}{t_{p + 1}}} } \int\limits_0^{{t_{p + 1}}} {\frac{{{{\ln }^m}\left( {1 + {t_{p + 2}}} \right)}}{{{t_{p + 2}}}}d{t_{p + 2}}} \nonumber \\
& = {\left( { - 1} \right)^{m + k + p}}m!k!\zeta \left( {\bar 2,{{\left\{ 1 \right\}}_{m - 1}}} \right)\zeta \left( {\bar 1,{{\left\{ 1 \right\}}_{p - 1}},2,{{\left\{ 1 \right\}}_{k - 1}}} \right)\nonumber \\
&\quad + {\left( { - 1} \right)^{m + k}}k!m!\zeta \left( {\bar 2,{{\left\{ 1 \right\}}_{k - 1}}} \right)\zeta \left( {\bar 1,{{\left\{ 1 \right\}}_{p - 1}},2,{{\left\{ 1 \right\}}_{m - 1}}} \right)\nonumber \\
&\quad + {\left( { - 1} \right)^{m + k + p}}k!m!\sum\limits_{i = 1}^{p - 1} {{{\left( { - 1} \right)}^i}\zeta \left( {\bar 1,{{\left\{ 1 \right\}}_{i - 1}},2,{{\left\{ 1 \right\}}_{m - 1}}} \right)\zeta \left( {\bar 1,{{\left\{ 1 \right\}}_{p - i - 1}},2,{{\left\{ 1 \right\}}_{k - 1}}} \right)} .
\end{align}
Hence, combining formulas (\ref{4.8})-(\ref{4.10}), then replacing $k-1$ by $k$, $m-1$ by $m$ and $p-1$ by $p$, we may easily deduce the desired result.\hfill$\square$\\
Setting $p=1,k=m=0$ in Theorem \ref{thm4.3}, we get
\[\zeta \left( {\bar 1,2,1,2} \right) + \zeta \left( {\bar 2,1,2} \right)\ln (2) + \zeta \left( {\bar 2} \right)\zeta \left( {\bar 1,1,2} \right) = \frac 1{2}\zeta^2\left({\bar 1},2 \right).\]
In fact, proceeding in a similar method to evaluation of the Theorems \ref{thm3.2}, \ref{thm4.1} and \ref{thm4.3}, it is possible to evaluate
other alternating multiple zeta values. For example, we have used our method to obtain the following explicit integral representations, closed form representations of multiple zeta values:
\begin{align}
&\zeta \left( {{{\left\{ {\bar 1} \right\}}_{2p + 2}},{{\left\{ 1 \right\}}_{k}}} \right) = {\left( { - 1} \right)^{p + 1}}{\rm{L}}{{\rm{i}}_{k + 2,{{\left\{ 2 \right\}}_p}}}\left( {\frac{1}{2}} \right)\nonumber\\
& = \frac{{{{\left( { - 1} \right)}^{k + p}}}}{{(k+1)!}}\int\limits_{0 < {t_{2p + 1}} <  \cdots  < {t_1} < 1} {\frac{{{{\ln }^{k+1}}\left( {1 - {t_{2p + 1}}} \right)}}{{\prod\limits_{j = 1}^p {\left\{ {\left( {1 + {t_{2j - 1}}} \right)\left( {1 - {t_{2j}}} \right)} \right\}\left( {1 + {t_{2p + 1}}} \right)} }}} d{t_1} \cdots d{t_{2p + 1}},\\
&\zeta \left( {{{\left\{ {\bar 1} \right\}}_{2p + 1}},{{\left\{ 1 \right\}}_{k}}} \right) = \frac{{{{\left( { - 1} \right)}^{k + p+1}}}}{{(k+1)!}}\int\limits_{0 < {t_{2p}} <  \cdots  < {t_1} < 1} {\frac{{{{\ln }^{k+1}}\left( {1 + {t_{2p}}} \right)}}{{\prod\limits_{j = 1}^p {\left\{ {\left( {1 + {t_{2j - 1}}} \right)\left( {1 - {t_{2j}}} \right)} \right\}} }}} d{t_1} \cdots d{t_{2p}},\\
&\sum\limits_{i = 0}^m {\frac{{{{\left( {\ln 2} \right)}^i}}}{{i!}}\zeta \left( {\bar 1,{{\left\{ 1 \right\}}_{m - i}},{{\left\{ {\bar 1} \right\}}_{2p}},{{\left\{ 1 \right\}}_{k}}} \right)} \nonumber \\
&= \frac{{{{\left( { - 1} \right)}^{k + p+1}}}}{{(k+1)!m!}}\int\limits_{0 < {t_{2p}} <  \cdots  < {t_1} < 0} {\frac{{{{\ln }^m}\left( {1 + {t_1}} \right){{\ln }^{k+1}}\left( {1 + {t_{2p}}} \right)}}{{\prod\limits_{j = 1}^p {\left\{ {\left( {1 + {t_{2j - 1}}} \right)\left( {1 - {t_{2j}}} \right)} \right\}} }}d{t_1} \cdots d{t_{2p}}}, \\
&\zeta \left( {\bar 1,{{\left\{ 1 \right\}}_m},{{\left\{ {\bar 1} \right\}}_{2p + 1}},{{\left\{ 1 \right\}}_{k}}} \right) = {\left( { - 1} \right)^{m + p + 1}}{\rm{L}}{{\rm{i}}_{k + 2,{{\left\{ 2 \right\}}_p},{{\left\{ 1 \right\}}_m}}}\left( {\frac{1}{2}} \right).
\end{align}
{\bf Acknowledgments.} We thank the anonymous referee for suggestions which led to improvements in the exposition.

 \end{document}